\documentclass[12pt]{article}

\usepackage{amssymb,amsmath,amsthm, amsfonts}

\def\sqr#1#2{{\vcenter{\vbox{\hrule height.#2pt
              \hbox{\vrule width.#2pt height#1pt \kern#1pt \vrule width.#2pt}
              \hrule height.#2pt}}}}
\def\signed #1{{\unskip\nobreak\hfil\penalty50
              \hskip2em\hbox{}\nobreak\hfil#1
              \parfillskip=0pt \finalhyphendemerits=0 \par}}
\def\endpf{\signed {$\sqr69$}}
\def\3n{\negthinspace \negthinspace \negthinspace }
\def\2n{\negthinspace \negthinspace }
\def\1n{\negthinspace }

\def\={\buildrel \triangle \over =}

\def\a{\alpha}
\def\b{\beta}

\def\th{\theta}

\def\G{\Gamma}
\def\D{\Delta}

\def\L{\Lambda}

\def\cA{{\cal A}}
\def\cB{{\cal B}}

\def\cF{{\cal F}}

\def\cL{{\cal L}}

\def\cY{{\cal Y}}

\def\ms{\medskip}

\def\hb{\hbox}

\def\esssup{\mathop{\rm esssup}}
\def\max{\mathop{\rm max}}
\def\min{\mathop{\rm min}}
\def\exp{\mathop{\rm exp}}
\def\sup{\mathop{\rm sup}}

\def\wt{\widetilde}
\def\cd{\cdot}

\def\esssup{\mathop{\hbox{\rm ess$\,$\rm sup$\,$}}}

\def\|{\Big |}
\def\({\Big (}
\def\){\Big )}
\def\[{\Big[}
\def\]{\Big]}
\def\be{\begin{equation}}
\def\bel{\begin{equation}\label}
\def\ee{\end{equation}}
\def\bt{\begin{theorem}}
\def\bcd{\begin{condition}}
\def\ecd{\end{condition}}
\def\et{\end{theorem}}
\def\bc{\begin{corollary}}
\def\ec{\end{corollary}}
\def\bde{\begin{definition}}
\def\ede{\end{definition}}
\def\bl{\begin{lemma}}
\def\el{\end{lemma}}
\def\bp{\begin{proposition}}
\def\ep{\end{proposition}}
\def\br{\begin{remark}}
\def\er{\end{remark}}
\def\ba{\begin{array}}
\def\ea{\end{array}}
\def\ed{\end{document}}

\def\square#1{\vbox{\hrule\hbox{\vrule height#1%
     \kern#1\vrule}\hrule}}
\def\rectangle#1#2{\vbox{\hrule\hbox{\vrule height#1%
     \kern#2\vrule}\hrule}}

\font\tenbb=msbm10 \font\sevenbb=msbm7 \font\fivebb=msbm5

\newfam\bbfam
\scriptscriptfont\bbfam=\fivebb \textfont\bbfam=\tenbb
\scriptfont\bbfam=\sevenbb

\newtheorem{lemma}{Lemma}[section]
\newtheorem{remark}{Remark}[section]

\newtheorem{theorem}{Theorem}[section]
\newtheorem{corollary}{Corollary}[section]

\newtheorem{definition}{Definition}[section]
\newtheorem{proposition}{Proposition}[section]
\newtheorem{condition}{Condition}[section]
\newtheorem{hypothesis}{Hypothesis}[section]

\makeatletter
   
   \@addtoreset{equation}{section}
\makeatother

\def\G{\Gamma}

\def\wt{\widetilde}
\def\wh{\widehat}
\def\essinf{\mathop{\mbox{\rm essinf }}}

\begin{document}

\title{Switching Game of Backward 
Stochastic \\ Differential Equations and Associated \\ System
of Obliquely Reflected Backward\\
Stochastic Differential Equations}

\author{Ying Hu\thanks{IRMAR, Universit\'e Rennes 1, Campus de Beaulieu,
35042 Rennes Cedex, France. Part of this work was completed when this author was visiting School of Mathematical Sciences,
Fudan University, whose hospitality is greatly appreciated. {\small\it E-mail:} {\small\tt Ying.Hu@univ-rennes1.fr}.\ms} \quad and \quad
Shanjian Tang\thanks{Corresponding author.  School of Mathematical
Sciences, Fudan University, Shanghai 200433, China.  Part of this work was completed when this author was visiting
IRMAR, Universit\'e Rennes 1, whose hospitality is greatly appreciated. {\small\it E-mail:} {\small\tt sjtang@fudan.edu.cn}.\ms}}

\maketitle

\abstract{This paper is concerned with the switching game of a one-dimensional backward
stochastic differential equation (BSDE). The associated Bellman-Isaacs equation is a system of matrix-valued BSDEs living in a special unbounded convex domain with  reflection on
the boundary along an oblique direction. In this paper, we show the existence of an adapted solution to this system of BSDEs with oblique reflection by the
penalization method, the monotone convergence, and the a priori estimates.}
\medskip


\medskip

{\bf Abbreviated title. }  Switching game of BSDEs and associated obliquely reflected
BSDEs



\newpage

\section{Introduction}

In this paper, we study the system of reflected BSDEs along an oblique direction arising naturally from the problem of switching game
of a scalar-valued BSDE. Let us first describe precisely the switching game problem by introducing some notations and hypotheses.

Let us fix a nonnegative real number $T>0$. First of all,
$W=\{W_t\}_{t\ge 0}$ is a standard Brownian motion with values in
$R^d$ defined on some complete probability space $(\Omega,{\cal
F},P)$. $\{{\cal F}_t\}_{t\ge 0}$ is the natural filtration of the
Brownian motion $W$ augmented by the $P$-null sets of ${\cal F}$.
All the measurability notions will refer to this filtration. In
particular, the sigma-field of predictable subsets of $[0,T]\times
\Omega$ is denoted by ${\cal P}$. Define $\Lambda:=\{1,\cdots,m_1\}$ and $\Pi:=\{1,\cdots,m_2\}$.

We denote by $S^2(R^{m_1\times m_2})$ or simply by $S^2$  the set of
$R^{m_1\times m_2}$-valued, adapted and c\`adl\`ag processes
$\{Y(t)\}_{t\in [0,T]}$ such that
$$||Y||_{S^2}:=E\left[\sup_{t\in [0,T]}|Y(t)|^2\right]^{1/2}<+\infty.$$
$(S^2,||\cdot||_{S^2})$ is then a Banach space.

We denote by $M^2((R^{m_1\times m_2})^d)$ or simply by $M^2$ the set
of (equivalent classes of) predictable processes $\{Z(t)\}_{t\in
[0,T]}$ with values in $(R^{m_1\times m_2})^d$ such that
$$||Z||_{M^2}:=E\left[\int_0^T |Z(s)|^2ds\right]^{1/2}<+\infty.$$
$M^2$ is then a Banach space endowed with this norm.

We define  also
\begin{eqnarray*}
N^2(R^{m_1\times m_2}):&=&\{K=(K_{ij})\in S^2 :\mbox{ for any }
(i,j)\in
\Lambda\times \Pi, K_{ij}(0)=0, \\
& &\mbox{ and } K_{ij}(\cdot) \mbox{ is increasing }\},
\end{eqnarray*}
which is abbreviated as $N^2$. $(N^2,||\cdot||_{S^2})$ is then a
Banach space.

Let $\psi$ be a
random function $\psi:[0,T]\times \Omega\times R\times R^d\times
\Lambda\times \Pi \to R$ whose component $\psi(\cdot,i,j)$ is
measurable with respect to ${\cal P}\otimes {\cal B}(R)\otimes
{\cal B}(R^d)$ for each pair $(i,j)\in \L\times \Pi$, and satisfies the following Lipschitz condition.
\begin{hypothesis}\label{psi}
(i) The generator $\psi(\cdot,0,0):=(\psi(\cdot,0,0,i,j))_{i\in
\L,j\in \Pi} \in M^2$.

(ii) There exists a constant $C> 0$ such that,  for each
$(t,y,y',z,z',i,j)\in [0,T]\times R\times R\times R^{d}\times
R^{d}\times\Lambda\times \Pi$,
$$|\psi(t,y,z,i,j)-\psi(t,y',z',i,j)|\le C (|y-y'|+|z-z'|), \quad a.s.$$
\end{hypothesis}

The  functions $k$
and $l$  are defined on $\Lambda\times \Lambda$ and
$\Pi\times \Pi$, respectively; their values are both positive. We make the following assumption on the  functions $k$ and $l$, which is standard in the literature.

\begin{hypothesis}\label{k}
(i) For $i\in \L$, $k(i,i)=0$.  For $(i,i')\in \Lambda\times
\Lambda$ such that $i\not=i'$, $k(i,i')> 0$.

(ii) For $j\in \Pi$, $l(j,j)=0$.  For $(j,j')\in \Pi\times \Pi$
such that $j\not=j'$, $l(j,j')> 0$.

(iii) For any $(i,i',i'')\in \Lambda\times \Lambda\times \Lambda$
such that $i\not=i'$ and $i'\not=i''$,
$$k(i,i')+k(i',i'')> k(i,i'').$$

(iv) For any $(j,j',j'')\in \Pi\times \Pi\times \Pi$
such that $j\not=j'$ and $j'\not=j''$,
$$l(j,j')+l(j',j'')> l(j,j'').$$

\end{hypothesis}

\bde  An admissible switching process for Player I ( resp. II ) on
$[t,T]$ with initial value $a_0\in \L$  ( resp. $b_0\in \Pi$ ) is defined to be
a pair of  sequences $\{a_i,\th_i\}_{i\ge 0}$ (resp.
$\{b_i,\tau_i\}_{i\ge 0}$  ), such that each $\th_i$ ( resp.
$\tau_i$ ) is an $\cF_{.}$-stopping time with
$$
t=\th_0\le \th_1\le \cdots\le T,\quad \hbox{\rm a.s.}
$$
$$\hbox{( resp. $ t=\tau_0\le \tau_1\le \cdots\le T,\quad \hbox{\rm
a.s.}$ ),}$$  each $a_i$ ( resp. $b_i$ ) is $\cF_{\th_i}$  (
resp. $\cF_{\tau_i}$ )  measurable with values in $\L$ ( resp.
$\Pi$ ),  and  there is an
integer-valued random variable $N(\cdot)$ satisfying
$$\theta_{N}=T \quad
\hbox{ \rm (resp. } \tau_N=T \hbox{\rm ) } \quad P\hbox{\rm -}
a.s. \quad \hbox{ \rm and } \quad N\in L^2(\cF_T).$$
  Denote by $\cA^a[t,\hat t]$ ( resp. $\cB^b[t, \hat t]$ )  the totality of the admissible switchings
for Player I  (resp. II)  on $[t, \hat t]$ with  the initial value
$a\in \L$ ( resp. $b\in \Pi$ ). Define the following abbreviations
for $t\in [0,T]$:
$$
{\cal A}^a_t:=\cA^a[t,T], \quad a \in \L; \quad
\cA_t:=\bigcup_{a\in \L}\cA_t^a
$$
and
$$
{\cal B}^b_t:=\cB^b[t,T], \quad b \in \Pi;\quad
\cB_t:=\bigcup_{b\in \Pi}\cB_t^b.
$$
 \ede

We shall identify $\{a_i,\th_i\}_{i\ge 0} \in \cA^a[t,T]$ with
\begin{equation}\label{admi}
a(s)=a_0\chi_{\{\theta_0\}}(s)+\sum_{i=1}^N a_{i-1}\chi_{(\th_{i-1},\th_i]}(s),\quad s\in
[t,T]. \end{equation}
 
For any $a(\cdot)\in {\cal A}_t$, we define the associated (cost)
process $A^{a(\cdot)}$ on $[t,T]$ as follows:
\begin{equation}\label{UA}
A^{a(\cdot)}(s)=\sum_{j=1}^{N-1}
k(a_{j-1},a_j)\chi_{[\theta_{j},T]}(s), \quad s\in
[t,T].
\end{equation}
Obviously, $A^{a(\cdot)}(\cdot)$ is a c\`adl\`ag process. In an
identical way, we define ${\cal B}^{b(\cdot)}$ on $[t,T]$ for
$b(\cdot)\in \cB_t$.

\bde  For $t\in [0,T]$ and $a\in \L ( \ \hb{\rm resp.}\  b\in \Pi\
)$, an admissible strategy $\a ^{a,t}$ ( resp. $\b^{b,t}$ ) with
the initial value $a\in \L$ (resp. $b\in \Pi$) for player $I$
(resp. $II$) on $[t,T]$ is a mapping $\a^{a,t}:\cup_{b\in
B}\cB^b[t,T]\to \cA^a[t,T]$ (resp. $\b^{b,t}:\cup_{a\in
A}\cA^a[t,T]\to \cB^b[t,T]$) such that
$$ b(s)=\wh b(s)\quad (\hb{\rm resp. } a(s)=\wh a(s)) \quad a.s. \ \forall s\in [t, \hat t],$$
implies
$$  \a^{a,t}[b(\cd)](s)=\a^{a,t}[\wh b(\cd)](s) \quad
(\hb{\rm resp. } \b^{b,t}[a(\cd)](s)=\b^{b,t}[\wh a(\cd)](s))
$$
for $s\in [t, \hat t]$.

We denote by $\G^a_t$ (resp. $\D^b_t$) all admissible strategies
with the initial value $a\in \L$ (resp. $b\in \Pi$) for player $I$
(resp. $II$) on $[t,T]$ . We adopt the convention that
$$
\cA^a_T=\{a\},\quad \G^a_T=\{a\}$$ and $$
        \cB^b_T=\{b\},\quad\D^b_T=\{b\}.
$$
\ede

Let $\xi$ be an $R^{m_1\times m_2}$-valued
${\cal F}_T$-measurable random variable. 
Now we are in position to introduce the switched BSDEs for both
players. For $t\in [0,T]$, $a(\cdot)\in {\cal A}_t$ and
$b(\cdot)\in {\cal B}_t$, consider the following BSDE: \be
\label{UV} \ba{rcl} U(s) &=&\displaystyle \xi_{a(T)b(T)}
+\left(A^{a(\cd)}(T)-A^{a(\cd)}(s)\right)-\left(B^{b(\cd)}(T)-B^{(\cd)}(s)\right)\\[3mm]
&&\displaystyle+\int_s^T\psi(r,U(r),V(r),a(r), b(r))\,
dr-\int_s^TV(r)\, dW(r), \quad s\in [t,T].\ea \ee This is a
(slightly) generalized BSDE: it is equivalent to the following
standard BSDE: \be \label{UV1}\ba{rcl} {\bar U}(s)
&=&\displaystyle\xi_{a(T)b(T)}
+A^{a(\cd)}(T)-B^{b(\cd)}(T)\\[3mm]
&&\displaystyle+\int_s^T\psi(r,\bar{U}(r)-A^{a(\cdot)}(r)+B^{b(\cdot)}(r),\bar{V}(r),a(r))\,
dr\\[3mm]
&&\displaystyle-\int_s^T\bar{V}(r)\, dW(r),\quad s\in [t,T]\ea\ee
via the simple change of variable:
$$\bar{U}(s)=U(s)+A^{a(\cdot)}(s)-B^{b(\cdot)}(s),\quad \bar{V}(s)=V(s).$$

Hence, for each pair $(a(\cdot), b(\cdot))\in \cA_t\times \cB_t$,
BSDE (\ref{UV}) has a unique solution in $S^2\times M^2$,
 which will be denoted by $\left(U^{a(\cdot),
b(\cdot)},V^{a(\cdot), b(\cdot)}\right)$. We note that $U$ is only
a c\`adl\`ag process.

The upper and lower switching game problems with the initial
scheme $(i,j)\in\Lambda\times \Pi$ are defined as follows:
$$
\mathop{\hbox{ \rm ess sup}}_{\b \in \D_t^j}\mathop{\hbox{ \rm ess
inf }}_{a(\cdot)\in \cA_t^i} U^{a(\cdot), \b(a(\cdot))}(t)
$$
and
$$
\mathop{\hbox{ \rm ess inf}}_{\a\in \G_t^i} \mathop{\hbox{ \rm ess
sup }}_{b(\cdot)\in \cB_t^j} U^{\a(b(\cdot)), b(\cdot)}(t),
$$
respectively. If
$$
\cY^{ij}(t):=\mathop{\hbox{ \rm ess sup}}_{\b \in
\D_t^j}\mathop{\hbox{ \rm ess inf }}_{a(\cdot)\in \cA_t^i}
U^{a(\cdot), \b(a(\cdot))}(t)=\mathop{\hbox{ \rm ess inf}}_{\a\in
\G_t^i} \mathop{\hbox{ \rm ess sup }}_{b(\cdot)\in \cB_t^j}
U^{\a(b(\cdot)), b(\cdot)}(t)
$$
for some $(i,j)\in \L\times \Pi$, we say that the switching game
with the initial scheme $(i,j)\in\Lambda\times \Pi$ has a value
$\cY^{ij}(t)$.

The above switching game (see, e.g. \cite{TangHou}) is associated to the 
following Bellman-Isaacs equation, which is a new type of reflected backward stochastic differential
equation (RBSDE for short)  with oblique reflection: for $(i,j)\in
\Lambda\times\Pi$  and $t\in [0,T]$,
\be\label{RBSDEi}\left\{\ba{rcl}Y_{ij}(t)&=&\displaystyle\xi_{ij}+\int_t^T
\psi(s,Y_{ij}(s),Z_{ij}(s),i,j)\,
ds\\[0.3cm]
& &\displaystyle-\int_t^TdK_{ij}(s)+\int_t^TdL_{ij}(s)-\int_t^TZ_{ij}(s)\, dW(s),\\[0.3cm]
Y_{ij}(t)&\le&\displaystyle \min_{i'\not=i}\{Y_{i'j}(t)+k(i,i')\}, \\[0.3cm]
Y_{ij}(t)&\ge&\displaystyle
\max_{j'\not=j}\{Y_{ij'}(t)-l(j,j')\},\\[0.3cm]
 &&\displaystyle
\int_0^T\left(Y_{ij}(s)-\min_{i'\not=i}\{Y_{i'j}(s)+k(i,i')\}\right)dK_{ij}(s)=0,\\[0.3cm]
& &\displaystyle\int_0^T\left(Y_{ij}(s)-\max_{j'\not=j}\{Y_{ij'}(s)-l(j,j')\}\right)dL_{ij}(s)=0. \ea\right.\ee Here,  the unknowns are the
processes $\{Y(t)\}_{t\in [0,T]}$, $\{Z(t)\}_{t\in [0,T]}$, $\{K(t)\}_{t\in [0,T]}$, and $\{L(t)\}_{t\in [0,T]}$, which are required to be
adapted with respect to the natural completed filtration of the Brownian motion $W$. Moreover, $K$ and $L$ are componentwisely increasing
processes. The last two relations in (\ref{RBSDEi}) are called the  upper and lower minimal boundary conditions.

One-dimensional RBSDEs were first studied by El Karoui et al.
\cite{ElKPPQ} in the
 case of one obstacle, and then by Cvitanic and Karatzas~\cite{CvtanicKaratzas}
 in the case of two obstacles. In both papers, it is recognized that one-dimensional reflected
 BSDEs, with one obstacle and with two obstacles, are generalizations of
optimal stopping and Dynkin games, respectively. Nowadays, the
literature on one-dimensional reflected BSDEs is very rich. The
reader is referred to  Peng
and Xu~\cite{PengXu} and Buckdahn and Li~\cite{bl}, among others, for the one-dimensional
reflected BSDEs with two obstacles.

 Multi-dimensional RBSDEs were studied by
Gegout-Petit and Pardoux \cite{Gegout-PetitPardoux}, but their
BSDE is reflected on the boundary of a convex domain along the
inward normal direction, and their method depends heavily on the
properties of this inward normal reflection (see (1)-(3) in
\cite{Gegout-PetitPardoux}). We note that in a very special case
(e.g., $\psi$ is independent of $z$), Ramasubramanian
\cite{Ramasubramanian} studied a BSDE in an orthant with oblique
reflection. Multi-dimensional BSDEs reflected along an oblique
direction rather than a normal direction, still remains to be open
in general, even in a convex domain, let alone in a nonconvex
domain. Note  that there are some papers dealing with SDEs with
oblique reflection (see, e.g. \cite{LionsSznitman,DupuisIshii}).

In our previous work~\cite{HuTang}, we studied the
optimal switching problem for \break one-dimensional BSDEs, and
the associated following type of obliquely reflected \break
multi-dimensional BSDEs: for $i\in \Lambda$,
\be\label{RBSDEswitching}\left\{\ba{rcl}Y_i(t)&=&\displaystyle\xi_i+\int_t^T
\psi(s,Y_i(s),Z_i(s),i)\,
ds\\[0.3cm]
& &\displaystyle-\int_t^TdK_i(s)-\int_t^TZ_i(s)\, dW(s),\\[0.3cm]
Y_i(t)&\le&\displaystyle \min_{i'\not=i}\{Y_{i'}(t)+k(i,i')\}, \\[0.3cm]
 &&\displaystyle
\int_0^T\left(Y_i(s)-\min_{i'\not=i}\{Y_{i'}(s)+k(i,i')\}\right)dK_i(s)=0.\\
\ea\right.\ee It should be added that a less
general form of RBSDE (\ref{RBSDEswitching}) (where  the generator
$\psi$ does not depend on $(y,z)$) is suggested by \cite{CarLud}.
But they did not discuss the existence and uniqueness of solution,
which is considered to be difficult. See Remark 3 in
\cite{CarLud}.


Recently, Tang, Zhong and Koo~\cite{TangZhongKoo} discussed the mixed
switching and stopping problem for one-dimensional BSDEs, and
obtained the existence and uniqueness result for  the
associated following type of multi-dimensional obliquely reflected
BSDEs: for $i\in \Lambda$ and $t\in [0,T]$,
\be\label{RBSDEswitchingStopping}\left\{\ba{rcl}Y_i(t)&=&\displaystyle\xi_i+\int_t^T
\psi(s,Y_i(s),Z_i(s),i)\,
ds\\[0.3cm]
& &\displaystyle-\int_t^TdK_i(s)+\int_t^T\, dL_i(s)-\int_t^TZ_i(s)\, dW(s),\\[0.3cm]
Y_i(t)&\le&\displaystyle \min_{i'\not=i}\{Y_{i'}(t)+k(i,i')\}, \quad Y_i(t)\ge S(t),\\[0.3cm]
 &&\displaystyle
\int_0^T\left(Y_i(s)-\min_{i'\not=i}\{Y_{i'}(s)+k(i,i')\}\right)dK_i(s)=0,\\[0.3cm]
&&\displaystyle \int_0^T(Y_i(t)-S(t))\, dL_i(t)=0.  \ea\right.\ee Here, $S$ is a  given $\{\cF_t, 0\le t\le T\}$-adapted process with some
suitable regularity.

RBSDE (\ref{RBSDEi})  is more complicated than that of RBSDE
(\ref{RBSDEswitching}) arising from the optimal switching
problem for BSDEs. For each fixed $j\in \Pi$, if we do not impose
the following constraint: \be Y_{ij}(t)\ge
\max_{j'\not=j}\{Y_{ij'}(t)-l(j,j')\}, \quad t\in [0,T],\ee and
its related boundary condition: \be
\int_0^T\left(Y_{ij}(s)-\max_{j'\not=j}\{Y_{ij'}(s)-l(j,j')\}\right)dL_{ij}(s)=0,\ee
then we can take $L\equiv 0$, and RBSDE (\ref{RBSDEi}) is reduced
to RBSDE (\ref{RBSDEswitching}).

 RBSDE (\ref{RBSDEi}) evolves in the closure $\overline{Q}$ of domain $Q$:
\be \ba{rcl} Q&:=&\displaystyle \biggl\{(y_{ij})\in
\mathrm{R}^{m_1\times m_2}: y_{ij}< y_{i'j}+k(i,i') \\
&&\displaystyle\qquad \mbox{\rm for any } i,i'\in \Lambda \mbox{
\rm such that }i'\not=i \hbox{
\rm and }  j\in \Pi; \\
&& \qquad \qquad \qquad\qquad \displaystyle y_{ij}>y_{ij'}-l(j,j') \\
&&\displaystyle\qquad  \mbox{\rm for any } j,j'\in \Pi \mbox{ \rm
such that }j'\not=j \hbox{ \rm and }  i\in \Lambda\biggr\},\ea\ee
which is convex and unbounded. The boundary $\partial Q$ of domain
$Q$ consists of the boundaries $\partial D_{ij}^-$ and $\partial
D_{ij}^+$, $(i,j)\in\Lambda\times \Pi,$ with
$$D_{ij}^-:=\{(y_{ij})\in \mathrm{R}^{m_1\times m_2}:y_{ij}<y_{i'j}+k(i,i'),
 \mbox{ for any }i'\in\Lambda \mbox{ such that } i'\not=i\}$$
 and
$$D_{ij}^+:=\{(y_{ij})\in \mathrm{R}^{m_1\times m_2}:y_{ij}>y_{ij'}-l(j,j'),
 \mbox{ for any }j'\in\Pi \mbox{ such that } j'\not=j\}$$
 for $(i,j)\in\Lambda\times \Pi.$
That is,
$$\partial Q=\mathop{\cup}_{i=1}^{m_1}\mathop{\cup}_{j=1}^{m_2}\left(\partial D_{ij}^-\cup \partial D_{ij}^+\right). $$
In the interior of $\overline {Q}$, each equation in
(\ref{RBSDEi}) is independent of others. On the boundary, say
$\partial D_{ij}^-$ (resp. $\partial D_{ij}^+$), the $(i,j)$-th
equation  is switched to another one $(i',j)$ (resp. $(i,j')$),
and the solution is reflected along the oblique direction
$-e_{ij}$ (resp. $e_{ij}$), which is the negative (resp. positive)
direction of the $(i,j)$-th coordinate axis.

The existence of solution for RBSDE (\ref{RBSDEi})
constitutes a main contribution of this paper. We prove the
existence by a penalization method.  Proving the existence of solution of RBSDE (\ref{RBSDEi}) presents new difficulties when
one follows our previous work~\cite{HuTang} using the penalization
method. In fact, in order to establish the a priori estimates which are essential for the proof of the existence, we have to
use the representation of solutions to obliquely reflected BSDEs proved in
\cite{HuTang}, and  we have to impose the additional technical condition that the generator $\psi$ is uniformly bounded.
The question of uniqueness is still an open problem.

The rest of the paper is organized as follows: in Section 2, we prove the existence of solution for RBSDE (\ref{RBSDEi})
by a penalization method.  The last section is devoted to discussions on some possible extensions.

\section{Existence of an adapted solution to the associated RBSDE}

In this section, we state and prove our existence result for RBSDE
(\ref{RBSDEi}).

We need the following additional technical assumption.

\begin{hypothesis}\label{UB} The generator $\psi$ is
uniformly bounded with respect to all its arguments.
\end{hypothesis}

we shall use $\|\psi\|_\infty$ to denote the least upper bound of
$|\psi|$.

\bde  An adapted solution to RBSDE  (\ref{RBSDEi}) is defined to
be a set $(Y, Z, K, L)=\{Y(t),Z(t),K(t), L(t)\}_{t\in [0,T]}$ of
predictable processes with values in $ \left(R^{m_1\times
m_2}\right)^{1+d+1+1}$ such that $P$-a.s., $t\mapsto Y(t)$ is continuous,
$t\mapsto K(t)$ and $t\mapsto L(t)$ are continuous and componentwisely increasing, $t\mapsto
Z(t)$ belongs to $L^2(0,T; (R^{m_1\times m_2})^d)$, $t\mapsto
\psi(t,Y_{ij}(t),Z_{ij}(t),i,j)$ belongs to $L^1(0,T; R^{m_1\times m_2})$
and $P$-a.s., RBSDE (\ref{RBSDEi}) holds for each $t\in [0,T]$.
\ede

The main result of this paper is  the following
existence of an adapted solution to RBSDE (\ref{RBSDEi}).

\begin{theorem}\label{existence}
Let Hypotheses \ref{psi}, \ref{k} and \ref{UB} be satisfied.
Assume that
$$\xi\in L^2(\Omega,{\cal F}_T,P;R^{m_1\times m_2})$$ takes values in
$\bar{Q}$. Then RBSDE  (\ref{RBSDEi}) has an adapted solution
$(Y,Z,K,L)$ in $S^2\times M^2\times (N^2)^2$.
\end{theorem}

We first sketch the proof.

{\bf Sketch of the Proof:} The proof is divided into five subsections.
In Subsection 3.1, we introduce the penalized RBSDEs whose existence of solution follows from a slightly generalized result in \cite{HuTang}. In Subsection 3.2, we give the (implicit) representation of these solutions. In Subsection 3.3, we state a fundamental lemma and some (uniform) a priori estimates for these
solutions. In Subsection 3.4, we prove the (monotone) convergence of these solutions. And the last subsection is devoted to checking out the boundary conditions.


\subsection{The penalized RBSDEs}

We shall use a penalization method to construct a solution to
RBSDE  (\ref{RBSDEi}). We observe (as mentioned in the
introduction) that RBSDE  (\ref{RBSDEi}) consists of the $m_2$
systems of $m_1$-dimensional obliquely reflected BSDEs of the form
like (\ref{RBSDEswitching}):
\be\left\{\ba{rcl}Y_{ij}(t)&=&\displaystyle\xi_{ij}+\int_t^T
\psi(s,Y_{ij}(s),Z_{ij}(s),i,j)\,
ds\\[0.3cm]
& &\displaystyle-\int_t^TdK_{ij}(s)+\int_t^TdL_{ij}(s)-\int_t^TZ_{ij}(s)\, dW(s),\\[0.3cm]
Y_{ij}(t)&\le&\displaystyle \min_{i'\not=i}\{Y_{i'j}(t)+k(i,i')\}, \\[0.3cm]
 &&\displaystyle
\int_0^T\left(Y_{ij}(s)-\min_{i'\not=i}\{Y_{i'j}(s)+k(i,i')\}\right)dK_{ij}(s)=0;
\quad i\in \L, \ea\right.\ee with the unknown processes being
$$(Y_{ij}, Z_{ij},K_{ij}; i=1,2, \ldots, m_1)$$
(the process $(L_{1j}, \ldots, L_{m_1j})$ is taken to be
previously given) for $j=1,2,\ldots,m_2$. These $m_2$ systems have
been well studied by Hu and Tang~\cite{HuTang}. In RBSDE
(\ref{RBSDEi}), they are coupled together by the processes
$(L_{1j}, \ldots, L_{m_1j})$ through the constraint \be
Y_{ij}(t)\ge \max_{j'\not=j}\{Y_{ij'}(t)-l(j,j')\}, \quad (i,j)\in
\L\times \Pi \ee and the minimal boundary condition: \be
\int_0^T\left(Y_{ij}(s)-\max_{j'\not=j}\{Y_{ij'}(s)-l(j,j')\}\right)dL_{ij}(s)=0,
\quad (i,j)\in \L\times \Pi. \ee Therefore, it is natural to
consider the following penalized system of RBSDEs (the unknown
processes are $(Y_{ij}, Z_{ij},K_{ij}; i\in \L, j\in \Pi)$:
\be\label{penalized}\left\{\ba{rcl}Y_{ij}(t)&=&\displaystyle\xi_{ij}+\int_t^T
\psi(s,Y_{ij}(s),Z_{ij}(s),i,j)\,
ds\\[0.3cm]
& &\displaystyle +n\sum_{j'=1}^{m_2}\int_t^T\left(Y_{ij}(s)-Y_{ij'}(s)+l(j,j')\right)^-\,
ds \\[0.3cm]
& &\displaystyle-\int_t^TdK_{ij}(s)-\int_t^TZ_{ij}(s)\, dW(s);\\
Y_{ij}(t)&\le&\displaystyle \min_{i'\not=i}\{Y_{i'j}(t)+k(i,i')\}; \\[0.3cm]
 &&\displaystyle
\int_0^T\left(Y_{ij}(s)-\min_{i'\not=i}\{Y_{i'j}(s)+k(i,i')\}\right)dK_{ij}(s)=0;
\quad (i,j)\in \L\times \Pi. \ea\right.\ee
Note that when $j'=j$, we have, in
view of Hypothesis \ref{k} (ii),
\be(Y_{ij}(s)-Y_{ij'}(s)+l(j,j'))^-=0. \ee



Note also that for any integer $n$, the $ij$-th component of the generator of
(\ref{penalized}) depends also on $y_{ij'}, j'\not= j$. Hence we cannot apply
directly the existence  result  in \cite{HuTang}. However,
by slightly adapting the relevant arguments in \cite{HuTang}, we have the following assertion.

\begin{proposition}\label{existence-penalized}  For any integer $n$, RBSDE (\ref{penalized}) has
an adapted solution $(Y^n,Z^n, K^n)$ in the space $S^2\times
M^2\times N^2$.
\end{proposition}
{\bf Proof.}
Let the integer $n$ be fixed. For an integer $m$, consider the following penalized BSDE whose solution is denoted by $(Y^{n,m},Z^{n,m})$:
\be\label{penalizednm}\left\{\ba{rcl}Y_{ij}(t)&=&\displaystyle\xi_{ij}+\int_t^T
\psi(s,Y_{ij}(s),Z_{ij}(s),i,j)\,
ds\\[0.3cm]
& &\displaystyle +n\sum_{j'=1}^{m_2}\int_t^T\left(Y_{ij}(s)-Y_{ij'}(s)+l(j,j')\right)^-\,
ds \\[0.3cm]
& &\displaystyle -m\sum_{i'=1}^{m_1}\int_t^T\left(Y_{ij}(s)-Y_{i'j}(s)-k(i,i')\right)^+\,
ds \\[0.3cm]
& &\displaystyle-\int_t^TZ_{ij}(s)\, dW(s); \quad (i,j)\in \L\times \Pi. \ea\right.\ee From the comparison theorem for multi-dimensional BSDEs
in \cite{HuPeng}, $\{Y_{ij}^{n,m}(t)\}_m$ is decreasing. Following the relevant arguments in \cite{HuTang}, we prove that there exists an
adapted solution $(Y^n,Z^n, K^n)$ in the space $S^2\times M^2\times N^2$, and moreover
$$Y_{ij}^n(t)=\lim_{m\rightarrow \infty} Y_{ij}^{n,m}(t).$$
\endpf

\subsection{Representation and uniqueness of  solution to the penalized RBSDE}

Note again that for any integer $n$, the $ij$-th component of the generator of
(\ref{penalized}) depends also on $y_{ij'}, j'\not= j$. Hence we cannot apply
directly the Representation Theorem 3.1  in \cite{HuTang}.

Nevertheless, by defining a new generator
\be
 \wt \psi (r, y, z, i, j;n):= \psi (r, y, z, i,
j)+n\sum_{j'\not=j}\left(y-Y_{ij'}^n(r)+l(j,j')\right)^-\ee for
$(r,y,z)\in [t,T]\times R\times R^d,$ and $(i,j)\in \L\times \Pi$ to include
$Y^n_{ij'},j'\not=j$ in it, we have for each $j\in
\Pi$, the triplet $(Y^n_{ij},Z^n_{ij}, K^n_{ij}; i\in \L)$ is an
adapted solution of the following $m_1$-dimensional RBSDE:
\be\left\{\ba{rcl}Y_{ij}(t)&=&\displaystyle\xi_{ij}+\int_t^T
\wt\psi(s,Y_{ij}(s),Z_{ij}(s),i,j;n)\,
ds\\[0.3cm]
& &\displaystyle-\int_t^TdK_{ij}(s)-\int_t^TZ_{ij}(s)\, dW(s);\\[0.3cm]
Y_{ij}(t)&\le&\displaystyle \min_{i'\not=i}\{Y_{i'j}(t)+k(i,i')\}; \\[0.3cm]
 &&\displaystyle
\int_0^T\left(Y_{ij}(s)-\min_{i'\not=i}\{Y_{i'j}(s)+k(i,i')\}\right)dK_{ij}(s)=0;
\quad i\in \L. \ea\right.\ee
Then we can apply the Representation Theorem 3.1  in \cite{HuTang}.

In order to state this representation theorem, first we introduce  some
notations.

Let $\{\theta_j\}_{j=0}^\infty$ be an increasing sequence of
stopping times with values in $[t,T]$ and $\forall j$, $\alpha_j$
is an ${\cal F}_{\theta_j}$-measurable random variable with values
in $\Lambda$, and $\chi$ is the indicator function. We define
$$
a(s):=\alpha_0\chi_{\{\theta_0\}}(s)+\sum_{j=1}^\infty
\alpha_{j-1}\chi_{(\theta_{j-1},\theta_j]}(s), \quad s\in [t,T].
$$
The sequence $\{\theta_j, \alpha_j\}_{j=0}^\infty$ or $a(\cdot)$
is said to be an admissible switching strategy starting from the
mode $\alpha_0$, if there exists an integer-valued random variable
$N$ such that $\theta_{N}=T$, $P$-a.s. and $N\in L^2(\cF_T)$.

We denote by ${\cal A}_t$ the set of all these admissible switching
strategies and by ${\cal A}^i_t$ the subset of ${\cal A}$ consisting
of admissible switching strategies starting from the mode $i$.

For any $a(\cdot)\in {\cal A}_t$, we define the associated (cost)
process $A^{a(\cdot)}$ as follows:
$$
A^{a(\cdot)}(s)=\sum_{j=1}^{N-1}
k(\alpha_{j-1},\alpha_j)\chi_{[\theta_{j},T]}(s),\ s\in [t,T].
$$
Obviously, $A^{a(\cdot)}(\cdot)$ is an adapted increasing c\`adl\`ag process, and
$A^{a(\cdot)}(T)\in L^2({\cal F}_T)$ thanks to the fact that $N\in L^2({\cal F}_T)$.

Now we are in position to introduce the switched BSDE.
For $(t,i,j)\in [0,T)\times \L\times \Pi$ and
$a(\cdot)\in \cA^i_t$, $(U_j^{a(\cdot),n},V_j^{a(\cdot),n})$ is
the unique solution to the following BSDE:
\be\label{uj1}\ba{rcl} U_j(s)&=&
\displaystyle\xi_{a(T)j}+[A^{a(\cdot)}(T)-A^{a(\cdot)}(s)]
+\int_s^T \wt\psi(r, U_j(r),V_j(r),a(r), j;n)\, dr\\
&&\displaystyle -\int_s^TV_j(r)\, dW(r), \quad s\in [t,T].\ea \ee

On the other hand, for $a(\cdot)\in \cA_t^i$ of the following
form:\be \label{a()} a(s)=i\chi_{t}(s)+\sum_{p=1}^N
\alpha_{p-1}\chi_{(\theta_{p-1}, \theta_p]}(s), \quad s\in
[t,T],\ee we define for $s\in [t,T]$ and $j\in \Pi$,
\begin{eqnarray}\label{yj1}
 \wt Y^{a(\cd),n}_{j}(s)&:=&\sum_{p=1}^N
Y^n_{\a_{p-1},j}(s)\chi_{[\th_{p-1},
\th_p)}(s)+\xi_{a(T)}\chi_{\{T\}}(s), \\
\wt Z^{a(\cd),n}_{j}(s)&:=&\sum_{p=1}^N
Z^n_{\a_{p-1},j}(s)\chi_{[\th_{p-1}, \th_p)}(s),\nonumber\\
\wt K^{a(\cdot),n}_{j}(s)&:=&\sum_{p=1}^N \int_{\th_{p-1}\wedge
s}^{\th_p\wedge s}d K^n_{\alpha_{p-1},j}(r),\nonumber
 \end{eqnarray}
and
\be\label{UAtilde}
\wt A^{a(\cdot),n}_{j}(s)=\sum_{p=1}^{N-1}[Y^n_{\alpha_p,
j}(\theta_p)+k(\alpha_{p-1},\alpha_p)- Y^n_{\alpha_{p-1},
j}(\theta_p)]\chi_{[\theta_p,T]}(s), \quad s\in [t,T].\ee
$\wt A^{a(\cdot),n}_{j}$
is an increasing process due to the fact that
$Y^n(s)$ satisfies the boundary condition in (\ref{penalized}), $\forall s\in [t,T].$ Then, for each $j\in \Pi$,  the triplet $(\wt Y^{a(\cd),n}_{j},\wt
Z^{a(\cd),n}_{j}, \wt K_{j}^{a(\cdot),n})$ is a solution to the
following BSDE:
\begin{eqnarray}\label{yj2}
& &\wt Y^{a(\cd),n}_{j}(s)\\
&=&\xi_{a(T) j}-[(\wt K^{a(\cdot),n}_{j}(T)+\wt
A^{a(\cdot),n}_{j}(T))- (\wt K^{a(\cdot),n}_{j}(s)+\wt
A^{a(\cdot),n}_{j}(s))]
+A^{a(\cdot)}(T)-A^{a(\cdot)}(s)\nonumber\\
&&\displaystyle+\int_s^T\wt\psi(r,\wt Y^{a(\cd),n}_{j}(r),\wt
Z^{a(\cd),n}_{j}(r),a(r),j;n)\, dr-\int_s^T\wt
Z^{a(\cd),n}_{j}(r)\, dW(r). \nonumber
\end{eqnarray}

 As
$$Y^n_{a(s)j'}(s)=\wt{Y}_{j'}^{a(\cdot),n}(s), \quad a.e. \quad
s\in [t,T],$$ BSDE~(\ref{uj1}) can be rewritten as the
following equation:
 \be\label{uj2}\ba{rcl} U_j(s)&=&
\displaystyle\xi_{a(T)j}+[A^{a(\cdot)}(T)-A^{a(\cdot)}(s)]
+\int_s^T \psi(r, U_j(r),V_j(r),a(r), j)\, dr\\
&&\displaystyle
+n\sum_{j'\not=j}\int_s^T(U_j(r)-\wt{Y}_{j'}^{a(\cdot),n}(r)+l(j,j'))^-\,
dr-\int_s^TV_j(r)\, dW(r),\\
&&\displaystyle \qquad\qquad \quad s\in [t,T].\ea \ee

We are now ready to state the representation formula which is taken from
Theorem 3.1 in \cite{HuTang}.

\begin{proposition} \label{representation} Assume that $a(\cdot)\in \cA_t^i$. Then we have
\be \label{iner3.2}\esssup_{a(\cdot)\in \cA_t^i}\left(\wt Y^{a(\cdot),n}_j(s)-
U^{a(\cdot),n}_{j}(s)\right)=0,\quad s\in [t,T], \quad j\in \Pi, \ee
 which, putting in particular  $s=t$, implies
$$Y^n_{ij}(t)=\essinf_{a(\cdot)\in \cA_t^i}U^{a(\cdot),n}_{j}(t), \quad j\in \Pi.$$
\end{proposition}

It is crucial to observe  that the above representation formula is implicit since $U^{a(\cdot),n}_{j}$ still depends on $Y^n$. However, it is
sufficient for us to deduce the a priori estimates. Also, it is sufficient for us to deduce the uniqueness of the solution to the penalized
RBSDE~(\ref{penalized}). In fact, if we have two solutions $(Y^{k,n}_{ij},Z^{k,n}_{ij}, K^{k,n}_{ij}; i\in \L)$ with $k=1,2$, we can define
$(U_j^{k, a(\cdot),n},V_j^{k, a(\cdot),n}), k=1,2,$ for $(t,i,j)\in [0,T)\times \L\times \Pi$ and $a(\cdot)\in \cA^i_t$, as the unique
solution to BSDE~(\ref{uj1}) with $Y^{n}_{ij}$ being replaced with $Y^{1,n}_{ij}$ and $Y^{2,n}_{ij}$, respectively. We have the estimate for
$U_j^{1, a(\cdot), n}-U_j^{2, a(\cdot), n}$ :
$$|U_j^{1, a(\cdot), n}(t)-U_j^{2, a(\cdot), n}(t)|^2\le C_n E[\int_t^T |Y^{1,n}(s)-Y^{2,n}(s)|^2 ds|F_t]$$
for some constant $C_n$. Then applying  the above representation
and Gronwall inequality, we obtain that $Y^{1,n}=Y^{2,n}$, which
gives the uniqueness.

\subsection{A basic lemma and a priori estimates}

 For each integer $n$, let $(Y^n,Z^n, K^n)\in S^2\times M^2\times N^2$ be the  adapted
solution of RBSDE (\ref{penalized}).
Intuitively, as $n$ tends to $+\infty$, we expect that the
sequence of solutions
$$\{(Y^n,Z^n, K^n)\}_{n=1}^\infty$$ together with the penalty term
$$
L^n_{ij}(t):=n\sum_{j'=1}^{m_2}\int_0^t\left(Y_{ij}(s)-Y_{ij'}(s)+l(j,j')\right)^-\,
ds, \quad (t,i,j)\in [0,T]\times \L\times \Pi
$$
will have a limit $(Y,Z, K, L)$, which solves RBSDE
(\ref{RBSDEi}).

For this purpose, it is crucial to prove that the  penalty term is
bounded in some suitable sense. Then we are naturally led to
compute
$$ \left(Y_{ij}(t)-Y_{ij'}(t)+l(j,j')\right)^-,$$
using It\^o-Meyer's formula, as done in  \cite{HuTang}. However, in our present situation, the
additional term $K^n$ appears in RBSDE (\ref{penalized}), which
gives rise to a serious difficulty  to derive the bound of $L^n$
in the preceding procedure. In what follows, we shall use the
representation result for $Y^n$  of Proposition~\ref{representation} to get
around the difficulty.

We have the following basic lemma.

\bl\label{basic} For $j,j'\in\Pi$ and $a(\cdot)\in \cA_t^i$,  we have
\be\label{key}
n\left(U_{j}^{a(\cdot),n}(s)-\wt{Y}_{j'}^{a(\cdot),n}(s)+l(j,j')\right)^-\le
2\|\psi\|_\infty, \quad s\in [t,T].\ee Here, $U_{j}^{a(\cdot),n}$
and $\wt{Y}_{j'}^{a(\cdot),n}$ are defined by (\ref{uj1}) and (\ref{yj1}), respectively.
 \el

{\bf Proof. } We suppress the superscripts
 ${(a(\cdot),n)}$ of $U_{j}^{a(\cdot),n}, U_{j'}^{a(\cdot),n}, V_{j}^{a(\cdot),n}, V_{j'}^{a(\cdot),n}$, $\wt Y_{j}^{a(\cdot),n}$, $\wt Y_{j'}^{a(\cdot),n}$, $ \wt Y_{j''}^{a(\cdot),n}$
 and $\wt Z_{j'}^{a(\cdot),n}$ for simplicity.
 The whole proof consists of the following two steps.

\medskip
 {\bf Step 1. Calculation of the process $\overline
{U}_{jj'}(s)^-$, where $\overline
{U}_{jj'}(s):=U_{j}(s)-{\wt Y}_{j'}(s)+l(j,j'), s\in [t,T]$ using
It\^o-Meyer's formula.}

\medskip
In view of (\ref{uj2}) and (\ref{yj2}),  the process $\overline {U}_{jj'}(s),
s\in [t,T]$ satisfies the
following BSDE: \be\ba{rcl} &&\displaystyle \overline {U}_{jj'}(s)\\
&=&\displaystyle \overline {U}_{jj'}(T)+\int_s^T \left[\psi(r,U_{j}(r), V_{j}(r),a(r),j)-\psi(r,{\wt Y}_{j'}(r), {\wt Z}_{j'}(r), a(r),
j')\right]\, dr
\\[0.3cm]
&&\displaystyle+n\sum_{j''\not=j}\int_s^T (U_j-\wt{Y}_{j''}+l(j,j^{''}))(r)^-\, dr
-n\sum_{j''\not=j'}\int_s^T({\wt Y}_{j'}-\wt{Y}_{j''}+l(j',j''))(r)^-\, dr\\[0.6cm]
&&\displaystyle +\int_s^Td\left(\wt K^{a(\cdot),n}_{j}(r)+\wt A^{a(\cdot),n}_{j}(r)\right) -\int_s^T(V_{j}(r)-{\wt Z}_{j'}(r))\, dW(r), \quad
s\in [t,T]. \ea \ee

Applying It\^o-Meyer's formula (see, e.g. Meyer~\cite{Meyer}), we have
 \be\ba{rcl}
&&\displaystyle \overline {U}_{jj'}(s)^-+n\sum_{j''\not=j}\int_s^T
\chi_{\cL_{jj'}^-}(r)(U_j-\wt{Y}_{j''}+l(j,j^{''}))(r)^-\,
dr\\[0.3cm]
&&\displaystyle-n\sum_{j''\not=j'}\int_s^T\chi_{\cL_{jj'}^-}(r)({\wt
Y}_{j'}-\wt{Y}_{j''}+l(j',j''))(r)^-\,
dr+\int_s^T\, d{\widehat L}_{jj'}(r)\\[0.6cm]
&=&\displaystyle -\int_s^T\chi_{\cL_{jj'}^-}(r) \left[\psi(r,U_{j}(r), V_{j}(r),a(r),j)-\psi(r,{\wt Y}_{j'}(r), {\wt Z}_{j'}(r),
a(r),j')\right]\, dr
\\[0.3cm]
&&\displaystyle -\int_s^T\chi_{\cL_{jj'}^-}(r-)\, d\left(\wt K^{a(\cdot),n}_{j}(r)+\wt A^{a(\cdot),n}_{j}(r)\right)\\[0.3cm]
&&\displaystyle +\int_s^T\chi_{\cL_{jj'}^-}(r-)(V_{j}(r)-{\wt Z}_{j'}(r))\, dW(r), \quad s\in [t,T]\ea \ee where \be
\ba{c}\cL_{jj'}^-:=\{(s,\omega)\in [t,T]\times \Omega:\overline {U}_{jj'}(s)< 0\}, \ea \ee and ${\widehat L}_{jj'}$ is a c\`adl\`ag increasing
process. The above equation can be rewritten as
  \be\label{Tanaka}\ba{rcl}
&&\displaystyle \overline {U}_{jj'}(s)^-+n\int_s^T
\chi_{\cL_{jj'}^-}(U_j-\wt{Y}_{j'}+l(j,j^{'}))(r)^-\,
dr+\int_s^T \, d{\widehat L}_{jj'}(r)\\[0.5cm]
&=&\displaystyle -\int_s^T
\chi_{\cL_{jj'}^-}(r)\left[\psi(r,U_{j}(r),
V_{j}(r),a(r),j)-\psi(r,{\wt Y}_{j'}(r), {\wt Z}_{j'}(r), a(r),
j')\right]\,
dr \\[0.5cm]
&&\displaystyle -\int_s^T\chi_{\cL_{jj'}^-}(r-)\, d(\wt K^{a(\cdot),n}_{j}(r)+\wt A^{a(\cdot),n}_{j}(r))\\[0.5cm]
&&\displaystyle+\int_s^T\chi_{\cL_{jj'}^-}(r-)\left(V_{j}(r)-{\wt
Z}_{j'}(r)\right)\,
dW(r)\\[0.5cm]
&&\displaystyle+n\int_s^TI_{1,jj'}(r) \, dr+n\sum_{j''\not=
j,j''\not=j'}\int_s^T I_{2,jj'j''}(r)\, dr
 \ea \ee
 where the two integrands $I_{1,jj'}$ and $I_{2,jj'j''}$ are defined as follows
 \be I_{1,jj'}(r):=\chi_{\cL_{jj'}^-}(r)({\wt Y}_{j'}-\wt{Y}_{j}+l(j',j))^-, \quad r\in [t,T] \ee
 and
 \be I_{2,jj'j''}(r):=\chi_{\cL_{jj'}^-}(r)[({\wt Y}_{j'}-\wt{Y}_{j''}+l(j',j^{''}))(r)^-
-(U_{j}-\wt{Y}_{j''}+l(j,j''))(r)^-], \quad r\in [t,T].\ee

\medskip
{\bf Step 2. The BSDE for the process
$\{(U_j(s)-\wt{Y}_{j'}(s)+l(j,j'))^-, s\in [t,T]\}$. }

\medskip
In view of
Proposition \ref{representation},
$$\wt{Y}_{j}\le U_{j},$$
we have
\be
I_{1,jj'}\le \chi_{\cL_{jj'}^-}({\wt Y}_{j'}-U_j+l(j',j))^-=0, \quad j,j'\in\Pi,\ee
thanks to the fact that
$$l(j,j')+l(j',j)>l(j,j)=0.$$
Hence,
$$I_{1,jj'}(r)=0.$$

 Now we can rewrite (\ref{Tanaka}) as the following
equation:
 \begin{eqnarray}\label{Tanaka1}
& &(U_j(s)-\wt{Y}_{j'}(s)+l(j,j'))^-\nonumber\\
&=&n\sum_{j''\not=
j,j''\not=j'}\int_s^TI_{2,jj'j''}(r)\, dr+\int_s^T I_{3,jj'}(r)\,
dr
\nonumber\\
& & -\int_s^T \, d{\widehat L}_{jj'}(r) -\int_s^T\chi_{\cL_{jj'}^-}(r-)\, d(\wt K^{a(\cdot),n}_{j}(r)+\wt A^{a(\cdot),n}_{j}(r))\nonumber\\
& &-n\int_s^T
(U_j(r)-\wt{Y}_{j'}(r)+l(j,j^{'}))^- \,dr\nonumber\\
& &+\int_s^T\chi_{\cL_{jj'}^-}(r-)\left(V_{j}(r)-{\wt
Z}_{j'}(r)\right)\, dW(r),
\end{eqnarray}
where \be I_{3,jj'}(r):=\chi_{\cL_{jj'}^-}(r)\left[\psi(r,{\wt Y}_{j'}(r), {\wt Z}_{j'}(r), a(r), j')-\psi(r,U_{j}(r),
V_{j}(r),a(r),j)\right]\ee for $ r\in [t,T].$

We now show that \be I_{2,jj'j''}\le 0. \ee In fact, for
$j,j',j''\in\Pi$, taking into consideration  the elementary inequality that $x_1^--x_2^-\le (x_1-x_2)^-$,
for any two real numbers $x_1$ and $x_2$, we have \be \ba{rcl}
&&I_{2,jj'j''}=\chi_{\cL_{jj'}^-}[({\wt
Y}_{j'}-\wt{Y}_{j''}+l(j',j^{''}))(r)^-
-(U_{j}-\wt{Y}_{j''}+l(j,j''))(r)^-]\\
 &\le& \chi_{\cL_{jj'}^-} \left({\wt Y}_{j'}-U_{j}+l(j',j'')-l(j,j'')\right)^-= 0.
 \ea\ee
The last equality holds in the last relations, since $$ \{y\in
R^m:y_j-y_{j'}+l(j,j')<0\}\cap\{y\in
R^m:y_{j'}-y_{j}+l(j',j'')-l(j,j'')<0\}=\emptyset,$$
thanks to Hypothesis
\ref{k} (iv), i.e.,
$$l(j,j')+l(j',j'')> l(j,j'').$$

Define for $r\in [t,T]$, \begin{eqnarray}
I_{jj'}(r)&:=&n\sum_{j''\not=
j,j''\not=j'}\int_t^rI_{2,jj'j''}(s)\, ds
\nonumber\\
& & -\int_t^r \, d{\widehat L}_{jj'}(s)
-\int_t^r\chi_{\cL_{jj'}^-}(s-)\, d(\wt K^{a(\cdot),n}_{j}(s)+\wt
A^{a(\cdot),n}_{j}(s)).\nonumber \end{eqnarray} Obviously, the
process $I_{jj'}(\cdot)$ is a c\`adl\`ag decreasing process for
any $(j,j')\in \Pi\times \Pi$. BSDE~(\ref{Tanaka1}) is finally written
as the following equation:
\begin{eqnarray}\label{Tanaka2}
& &(U_j(s)-\wt{Y}_{j'}(s)+l(j,j'))^-\nonumber\\
&=&\int_t^T\, d I_{jj'}(r)+\int_s^T I_{3,jj'}(r)\, dr -n\int_s^T
(U_j(r)-\wt{Y}_{j'}(r)+l(j,j^{'}))^- \,dr\nonumber\\
& &+\int_s^T\chi_{\cL_{jj'}^-}(r-)\left(V_{j}(r)-{\wt
Z}_{j'}(r)\right)\, dW(r), \quad s\in [t,T].
\end{eqnarray}
Then, we have the following formula \begin{eqnarray}
(U_j(s)-\wt{Y}_{j'}(s)+l(j,j'))^-&=&E\left[\int_s^TI_{3,jj'}(r)\exp{[-n(r-s)]}\,
dr \biggm|\cF_s\right]\nonumber\\
&&+E\left[\int_s^T\exp{[-n(r-s)]}\, dI_{jj'}(r)
\biggm|\cF_s\right]\nonumber\\
&\le& 2\|\psi\|_\infty \int_s^T\exp{[-n(r-s)]}\, dr
.\end{eqnarray} Therefore, we have
$$n\left(U_{j}^{a(\cdot),n}(s)-\wt{Y}_{j'}^{a(\cdot),n}(s)+l(j,j')\right)^-\le
2\|\psi\|_\infty, \quad s\in [t,T].$$ This ends the proof. \endpf

\medskip
Thanks to this basic lemma, we deduce easily the following a priori estimates.

\begin{proposition}\label{apriori} (i) The sequence $\{Y_{ij}^n(t)\}_{n=1}^\infty$ is
increasing. Moreover, \be  -E\left[|\xi|
|\cF_t\right]-|\psi|_\infty T \le Y_{ij}^n(t) \le E\left[|\xi|
|\cF_t\right]+3|\psi|_\infty T ; \quad
E\left[\sup_t|Y_{ij}^n(t)|^2\right]\le C, \ee
where $C>0$ is a constant.

(ii) We have \be n\left(Y_{ij}^n(t)-Y_{ij'}^n(t)+l(j,j')\right)^-\le
2 \|\psi\|_\infty. \ee
\end{proposition}
{\bf Proof. }  (i)  According to the comparison theorem
for multi-dimensional BSDEs in \cite{HuPeng},
$$Y_{ij}^{n,m}(t)\le Y_{ij}^{n+1,m}(t),$$
where $Y_{ij}^{n,m}$ is defined via (\ref{penalizednm}).
Hence, the sequence $\{Y_{ij}^n(t)\}_{n=1}^\infty$ is
increasing by taking the limit when $m$ tends to $\infty$.

We have the following two facts in view of (\ref{uj2}):

(1) $U_{j}^{a(\cdot),n}(s)\ge -E\left[|\xi|
|\cF_s\right]-|\psi|_\infty T. $

(2) Taking $\bar a(\cdot)\equiv i,$ we have, from Lemma \ref{basic}, \be U_{j}^{\bar
a(\cdot),n}(s)\le  E\left[|\xi| |\cF_s\right]+|\psi|_\infty
T+2|\psi|_\infty T.\ee In view of the representation formula in Proposition \ref{representation}, we conclude the proof.

(ii) Putting $s=t$ in (\ref{key}), we obtain
$$n\left(U_{j}^{a(\cdot),n}(t)-Y_{ij'}^n(t)+l(j,j')\right)\ge
-2\|\psi\|_\infty.$$
From Proposition \ref{representation}, we deduce that
$$n\left(Y_{ij}^n(t)-Y_{ij'}^n(t)+l(j,j')\right)\ge -
2 \|\psi\|_\infty,$$
and the proof is complete.
\endpf

\subsection{Convergence of solutions}

We first prove that $(Z^n,K^n)$ is bounded.

\bl The pair of processes $(Z_{ij}^n, K_{ij}^n)$ are uniformly bounded in $M^2\times N^2$ for $(i,j)\in \L\times \Pi$.
\el

{\bf Proof. } From the RBSDE for $Y_{ij}^n$, in view of
Hypothesis~\ref{UB}, using It\^o's formula and Proposition~\ref{apriori}, we have \be\ba{rcl} &&\displaystyle E|Y_{ij}^n(0)|^2
+E\int_0^T|Z_{ij}^n(s)|^2\,
ds\\[0.5cm]
&\le &\displaystyle   E|\xi_{ij}|^2 +2\sum_{j'=1}^{m_2}
E\int_0^T|Y_{ij}^n(s)|
n\left(Y_{ij}^n(s)-Y_{ij'}^n(s)+l(j,j')\right)^-\, ds\\[0.5cm]
&&\displaystyle +
 2E\int_0^T|Y_{ij}^n(s)|\cdot
 |\psi(s,Y_{ij}^n(s), Z_{ij}^n(s),i,j)|\, ds \\[0.5cm]
&&\displaystyle +2E\int_0^T|Y_{ij}^n(s)|\, dK_{ij}^{n}(s)\\[0.5cm]
&\le &\displaystyle  C + C E\int_0^T  |Y_{ij}^n(s)|\, ds
+2E\left[\sup_t|Y_{ij}^n(t)|K_{ij}^{n}(T)\right]\\[0.5cm]
 &\le&\displaystyle  C_\epsilon +\epsilon
E[(K_{ij}^{n}(T))^2] \ea \ee and \be\ba{rcl}&&\displaystyle  E[(K_{ij}^{n}(T))^2]\\[0.5cm]
&\le&\displaystyle  C E|\xi_{ij}|^2+CE|Y_{ij}^n(0)|^2\\[0.5cm]
&&\displaystyle
+C\sum_{j'=1}^{m_2}E\int_0^T \left[n\left(Y_{ij}^n(s)-Y_{ij'}^n(s)+l(j,j')\right)^-\right]^2\, ds  \\[0.5cm]
&&\displaystyle +C E\int_0^T
  |\psi(s,Y_{ij}^n(s), Z_{ij}^n(s),i,j)|^2\, ds
+CE\int_0^T|Z_{ij}^n(s)|^2\, ds\\[0.5cm]
&\le &\displaystyle C +CE\int_0^T|Z_{ij}^n(s)|^2\, ds. \ea\ee
Combining the above two inequalities by taking a sufficiently small $\epsilon>0$, we conclude the proof.
\endpf

Define \be \b_{ij}^n(s):=
n\sum_{j'=1}^{m_2}\left(Y_{ij}^n(s)-Y_{ij'}^n(s)+l(j,j')\right)^-.\ee
Then \be L_{ij}^{n}(t)= \int_0^t\b_{ij}^n(s)\, ds. \ee

Next, we prove that $K^n$ is absolutely continuous whose derivative is uniformly bounded.

\bl  \label{5} For $(i,j)\in \L\times \Pi$ and an integer $n$, there is a uniformly bounded process  $\a_{ij}^n$  such that
$K_{ij}^{n}$ has the following form:
\be
K_{ij}^{n}(t)= \displaystyle  \int_0^t\a_{ij}^n(s)\, ds, \quad t\in [0,T].\ee
\el

{\bf Proof. } Fix the integer $n$. Consider the following penalized BSDEs:
 \be\ba{rcl} Y_{ij}(t) &=&\displaystyle
 \xi_{ij}+\int_t^T\left[\psi(s,Y_{ij}^n(s),Z_{ij}^n(s),i,j)+\b_{ij}^n(s)\right]\,
 ds
 \\[3mm]
 &&\displaystyle-m\sum_{i'=1}^{m_1}\int_t^T\left(Y_{ij}-Y_{i'j}-k(i,i')\right)^+\,
 ds\\[4mm]
 &&\displaystyle -\int_t^TZ_{ij}(s)\, dW(s),\ea \ee

 with $(i,j)\in \L\times \Pi$  and $t\in [0,T]$.
 It has a unique solution, denoted by $(\bar{Y}_{ij}^{n,m}, {\bar Z}_{ij}^{n,m})$.

 Proceeding similarly (in fact, much more simpler)
 as in Lemma~\ref{basic}, we can prove that for a constant $C>0$, \be
 \a^{n,m}_{ij}:=m\sum_{i'=1}^{m_1}\left(\bar{Y}_{ij}^{n,m}-\bar{Y}_{i'j}^{n,m}-k(i,i')\right)^+\le C. \ee
  Therefore, $\{\a_{ij}^{n,m}\}_{m=1}^\infty$   has a weak limit in
  $M^2$, denoted by $\a_{ij}^n$.  Then $\a_{ij}^n$
  is also uniformly bounded by the same constant $C$.

Define
\be
 \displaystyle \bar{K}_{ij}^{n,m}(t):= \int_0^t
 \a_{ij}^{n,m}(s)\, ds, \quad t\in [0,T].\ee
 From~\cite{HuTang}, we have \be\ba{rcl}\displaystyle
 \lim_{m\to
 \infty}\bar{Y}_{ij}^{n,m}(t)&=&Y_{ij}^{n}(t), \\
  \displaystyle \lim_{m\to
 \infty}\bar{Z}_{ij}^{n,m}(t)&=&Z_{ij}^{n}(t), \\
 \displaystyle \lim_{m\to
\infty}\bar{K}_{ij}^{n,m}(t)&=&\displaystyle \int_0^t\a_{ij}^n(s)\, ds=K_{ij}^{n}(t). \ea \ee
\endpf

Now we are able to state the convergence result. \bl    The sequence $\{Y^n, Z^n\}$ has a strong limit $(Y,Z)$  in $S^2\times M^2$.
The two sequences $\{\a^n\}$ and $\{\b^n\}$ have
 subsequences which converge to $\a$ and $\b$ weakly in $M^2$, respectively.\el

 {\bf Proof. } Note that $Y^n_{ij}$ is increasing in $n$. In
 view of Proposition~\ref{apriori} and applying the dominated convergence
 theorem, we deduce easily the strong convergence of  $\{Y^n\}$ in the
 space $M^2$. Note that  $(Y^n, Z^n)$  solves the following BSDE:
 \be\label{app}\ba{rcl}\displaystyle  Y_{ij}^n(t)&=&\displaystyle  \xi_{ij}
 +\int_t^T\left[\psi(s, Y_{ij}^n(s), Z_{ij}^n(s),i,j)+\b_{ij}^n(s)-\a_{ij}^n(s)\right]\,
ds\\[3mm]
 &&\displaystyle -\int_t^TZ_{ij}^n(s)\, dW(s),\quad t\in [0,T], \quad (i,j)\in \L\times \Pi, \ea \ee
 with $\{\a^n\}$ and $\{\b^n\}$  being uniformly bounded.

We now prove the strong convergence of $Z^n$.
Using It\^o's formula, we have
\be\ba{rcl}&&\displaystyle  |Y_{ij}^n(0)-Y_{ij}^m(0)|^2
+E\int_0^T|Z^n_{ij}(s)-Z^m_{ij}(s)|^2\, ds\\
&=&\displaystyle
2E\int_0^T\left(Y_{ij}^n(s)-Y_{ij}^m(s)\right)\left(\psi(s,Y_{ij}^n(s),Z_{ij}^n(s),
i,j)-\a_{ij}^n(s)+\b_{ij}^n(s)\right)\, ds\\
&&\displaystyle-
2E\int_0^T\left(Y_{ij}^n(s)-Y_{ij}^m(s)\right)\left(\psi(s,Y_{ij}^m(s),Z_{ij}^m(s),
i,j)-\a_{ij}^m(s)+\b_{ij}^m(s)\right)\, ds, \ea \ee
from which we deduce immediately that \be \lim_{n,m\to
\infty}E\int_0^T|Z^n(s)-Z^m(s)|^2\, ds=0. \ee

It
is routine to show the strong convergence of  $\{Y^n\}$ in the
space $S^2$.

Since $\{\alpha^n\}$ and $\{\beta^n\}$ are uniformly bounded,
 the last assertion of the lemma is obvious. \endpf

\medskip\bigskip

Define  for $(i,j,t)\in \L\times \Pi\times [0,T]$,
\be K_{ij}(t):=\int_0^t\a_{ij}(s)\, ds, \quad L_{ij}(t):=\int_0^t\b_{ij}(s)\, ds \ee
and
\be K:=(K_{ij}), \quad L:=(L_{ij}).\ee
Finally, we shall show that $(Y,Z,K,L)$  solves RBSDE (\ref{RBSDEi}).

In fact, it suffices to take the weak limit in $L^2(\cF_T)$ in BSDE~(\ref{app})
along a suitable subsequence to deduce that
 $(Y,Z,K,L)$ solves the following BSDE:
\be\ba{rcl}\displaystyle  Y_{ij}(t)&=&\displaystyle  \xi_{ij}
+\int_t^T\psi(s, Y_{ij}(s), Z_{ij}(s),i,j)\, ds-\int_t^T\,
dK_{ij}(s)\\[3mm]
&&\displaystyle +\int_t^T\, dL_{ij}(s)-\int_t^TZ_{ij}(s)\, dW(s),\quad (i,j)\in \L\times \Pi. \ea \ee

It remains to check out the boundary conditions, which will be given in the next subsection.

\subsection{Boundary conditions}

Let us first prove that $Y(t)\in \bar{Q}$.

On the one hand, as $(Y^n,Z^n,K^n)$ satisfies (\ref{penalized}),
we have
$$Y^n_{ij}(t)\le\displaystyle \min_{i'\not=i}\{Y^n_{i'j}(t)+k(i,i')\},$$
from which we deduce, by taking limit when $n$ tends to $\infty$, that
\begin{equation}\label{B1}Y_{ij}(t)\le\displaystyle \min_{i'\not=i}\{Y_{i'j}(t)+k(i,i')\}.\end{equation}

On the other hand, from Proposition \ref{apriori},
$$\left(Y_{ij}^n(t)-Y_{ij'}^n(t)+l(j,j')\right)^-\le
\frac{2 \|\psi\|_\infty}{n}.$$
Sending $n$ to $\infty$, we deduce that
\begin{equation}\label{B2}\left(Y_{ij}(t)-Y_{ij'}(t)+l(j,j')\right)^-=0.
\end{equation}

(\ref{B1}) and (\ref{B2}) shows that $Y(t)\in \bar{Q}$.

Now we check out the minimal boundary conditions.

From (\ref{penalized}), we have \be
E\int_0^T\left(Y_{ij}^n(s)-\min_{i'\not=i}\{Y_{i'j}^n(s)+k(i,i')\}\right)^-\a_{ij}^n(s)\,
ds=0.\ee Setting $n\to \infty$, we have \be
E\int_0^T\left(Y_{ij}(s)-\min_{i'\not=i}\{Y_{i'j}(s)+k(i,i')\}\right)^-\a_{ij}(s)\,
ds=0.\ee

On the other hand, from the construction, we have \be
E\int_0^T\left(Y_{ij}^n(s)-\max_{j'\not=j}\{Y_{ij'}^n(s)-l(j,j')\}\right)^+\b_{ij}^n(s)\,
ds=0.\ee Setting $n\to \infty$, we have \be
E\int_0^T\left(Y_{ij}(s)-\max_{j'\not=j}\{Y_{ij'}(s)-l(j,j')\}\right)^+\b_{ij}(s)\,
ds=0.\ee
The proof is then complete.

\section{Concluding remarks}

In this paper, we proved the existence of solution to RBSDE (\ref{RBSDEi}).
But the question of uniqueness remains open.

On the other hand, there exist different methods in the literature for the study of
switching control and game problems. For the classical method of
quasi-variational inequalities, the reader is referred to the book
of Bensoussan and Lions~\cite{BL}. See Tang and Yong~\cite{TY},
Pham, Ly Vath and Zhou~\cite{pz} and Tang and Hou~\cite{TangHou} and the
references therein for the theory of variational inequalities and
the dynamic programming for optimal stochastic switching control
and switching games. But these works are restricted to the
Markovian case. Recently, using the method of  Snell envelope
(see, e.g. El Karoui \cite{ElK}) combined with the theory of
scalar valued RBSDEs, Hamad\`ene and Jeanblanc \cite{HaJean}
studied the switching problem  in the non-Markovian context. The
obliquely reflected BSDE approach, first fully developed in Hu and
Tang~\cite{HuTang} for optimal stochastic switching and taking the
advantage of the theory and techniques of BSDEs, permits  to state
and solve these problems in a rather general non-Markovian
framework. The link between the solution of RBSDE (\ref{RBSDEi}) and the problem of switching games constitutes another challenge.


\begin{thebibliography}{99}

\bibitem{BL} A. Bensoussan and J. L. Lions,  Impulse Control
and Quasivariational Inequalities. {\it Gauthier-Villars,
Montrouge}, 1984.

\bibitem{bl} R. Buckdahn and J. Li, Stochastic differential games with reflection and related obstacle problems for Isaacs equations. {\it Acta Math. Appl. Sin. Engl. Ser.}  27 (2011), 647--678.

\bibitem{CarLud} R. Carmona and M. Ludkovski, Pricing asset scheduling flexibility using optimal switching.
{\it Appl. Math. Finance} 15 (2008), 405-447.


\bibitem{CvtanicKaratzas} J. Cvitanic and I. Karatzas, Backward
stochastic differential equations with reflection and Dynkin games. {\it Ann. Probab.
} 24 (1996), 2024--2056.


\bibitem{DupuisIshii} P. Dupuis and H. Ishii, SDEs with oblique reflection on nonsmooth domains. {\it
Ann.  Probab.} 21 (1993), 554--580.

\bibitem{ElK} N. El Karoui, Les aspects probabilistes du
contr\^ole stochastique. {\it Ninth Saint Flour Probability Summer
School -- 1979 (Saint Flour, 1979)}, pp. 73--238, Lecture Notes in
Math., 876, Springer, Berlin, 1981.

\bibitem{ElKPPQ} N. El Karoui, C. Kapoudjian, E. Pardoux, S. Peng
and M. C. Quenez, Reflected solutions of backward SDE's, and
related obstacle problems for PDE's. {\it Ann. Probab.} 25 (1997),
702--737.

\bibitem{Gegout-PetitPardoux}  A. Gegout-Petit and E. Pardoux,
Equations diff\'erentielles stochastiques r\'etrogrades
r\'efl\'echies dans un convexe. {\it Stochastics  Stochastic Rep.}
57 (1996), 111--128.


\bibitem{HaJean} S. Hamad\`ene and M. Jeanblanc, {\it On the starting
and stopping problem: application in reversible investments.}
Math. Oper. Res. 32 (2007), 182--192.

%

\bibitem{HuPeng} Y. Hu and S. Peng, On the comparison theorem
for multi-dimensional BSDEs. {\it C. R. Math. Acad. Sci. Paris}
343 (2006), 135--140.

\bibitem{HuTang} Y. Hu and S. Tang,  Multi-dimensional BSDE with oblique reflection and
optimal switching. {\it Probab. Theory  Related Fields} 147 (2010), 89--121.

\bibitem{LionsSznitman} P. L. Lions and A. S. Sznitman,  Stochastic differential equations
with reflecting boundary conditions. {\it Comm. Pure  Appl. Math.}
37 (1984), 511--537.

\bibitem{Meyer} P. A. Meyer, Un cours sur les int\'egrales
stochastiques. {\it S\'eminaire de Probabilit\'es, X}, pp.
245-400. Lecture Notes in Math., 511,  {\it Springer,
Berlin}, 1976.



\bibitem{PengXu} S. Peng and M. Xu, The smallest
$g$-supermartingale and reflected BSDE with single and double
$L^2$ obstacles. {\it Ann. Inst. H. Poincar\'e Probab. Statist.} 41 (2005),
605--630.

\bibitem{pz} H. Pham, V. Ly Vath and X. Y. Zhou, Optimal switching over multiple regimes. {\it SIAM J. Control Optim.} 48 (2009), 2217–-2253.

\bibitem{Ramasubramanian} S. Ramasubramanian,  Reflected backward stochastic differential equations in an orthant. {\it Proc.
Indian Acad. Sci. Math. Sci.} 112 (2002), 347--360.




\bibitem{TangHou} S. Tang and S. Hou,  Switching games of
stochastic differential systems, {\it SIAM J. Control
Optim.} 46 (2007), 900-929.

\bibitem{TY} S. Tang and J. Yong,  Finite horizon stochastic optimal
switching and impulse controls with a viscosity solution approach.
{\it Stochastics Stochastics Rep.} 45 (1993), 145--176.

\bibitem{TangZhongKoo} S. Tang, W. Zhong and H. Koo,  Optimal switching of one-dimensional reflected BSDEs and associated multidimensional BSDEs with oblique reflection, {\it SIAM J. Control Optim. } 49 (2011), 2279--2317.

\end{thebibliography}
\end{document}